\documentclass[12pt]{article}
\usepackage{amsfonts,amsmath,amssymb,amscd}
\usepackage[russian]{babel}
\usepackage{epsfig}
\long\def\comment#1\endcomment{}

\begin{document}

\begin{figure}[h]
\centerline{\uppercase{\bf Теоремы о галстуке. }}
\bigskip
\centerline{{\bf Крутовский Роман}
\footnote{rice12@yandex.ru }
}

{\bf Abstract.}

{\bf Theorem.} There are general position points $A, B, C, P$ on the projective plane. Let $A_P$ be the intersection point of lines $AP$ and $BC$. Analogously define $B_P$ and $C_P$. Take any points $A_1$, $B_1$, $C_1$ on $AP$, $BP$, $CP$, respectively. Let $W_C$ be the  intersection point of $A_1B_P$ and $B_1A_P$. Analogously define points $W_A$ and $W_B$. Then lines $CW_C$, $AW_A$ and $BW_B$ pass through one point.

We also generalize this theorem and find interesting related properties.

\bigskip
Следующая серия задач возникла в обсуждениях с Иваном Фроловым. Так же автор статьи очень благодарен Андрею Сергунину, сформулировавшему теорему $1.6$.

Отдельно хотелось бы поблагодарить Ф.\,А.\,Ивлева и А.\,Б.\,Скопенкова за поддержку и полезные замечания.

\smallskip
В статье используются базовые понятия и факты проективной геометрии. Подробную информацию о них можно найти в книгах [1] и [2].

{\it Соглашения.}

{\it 1)} Все теоремы рассматриваются в проективной плоскости.

{\it 2)} n-угольник --- n точек общего положения в плоскости (n --- натуральное число).

{\it 3)} Сторона n-угольника --- прямая, содержащая сторону n-угольника в обычном ее понимании.

{\it Примечание.} В тексте работы каждая из теорем, рисунок к ней и ее доказательство расположены друг за другом для удобства понимания самих теорем.

\smallskip
{\bf Определения.} Дан треугольник $ABC$. $P$ --- произвольная точка плоскости, не лежащая на сторонах треугольника $ABC$. $A_P$, $B_P$, $C_P$ --- точки пересечения прямых $AP$, $BP$ и $CP$ со сторонами  $BC$, $AC$, $AB$  соответственно. На $AP$, $BP$, $CP$ произвольным образом выбраны соответственно точки $A_1$, $B_1$, $C_1$. Через $W_C$ определим точку пересечения $A_1B_P$ и $B_1A_P$. Аналогично определяются точки $W_A$ и $W_B$. $W$ -- точка пересечения $AW_A$ и $BW_B$.

\smallskip
{\bf Теорема 1.1 (Теорема о галстуке).}  Прямая $CW_C$ проходит через $W$.

\smallskip
{\it Доказательство теоремы 1.1.} Для троек точек $A$, $A_1$, $A_P$ и $B$, $B_1$, $B_P$ применим теорему Паппа (см. рис. \ref{1}). Тогда получим, что $C$, $W_C$ и $J_C$, точка пересечения $AB_1$ и $BA_1$, лежат на одной прямой.  Аналогично определим точки $J_A$ и $J_B$ (и применим теорему Паппа для соответствующих им пар троек точек). По обратной теореме Брианшона конику можно вписать в $AC_1BA_1CB_1$ . Теперь заметим, что стороны шестиугольника $AC_1BA_1CB_1$ являются сторонами шестиугольника $AJ_CBJ_ACJ_B$, а значит, коника, вписанная в $AC_1BA_1CB_1$, вписана и в $AJ_CBJ_ACJ_B$. Тогда по теореме Брианшона $AJ_A, BJ_B, CJ_C$ пересекаются в некоторой точке~$W$.
QED
\end{figure}

\begin{figure}[h]

\includegraphics{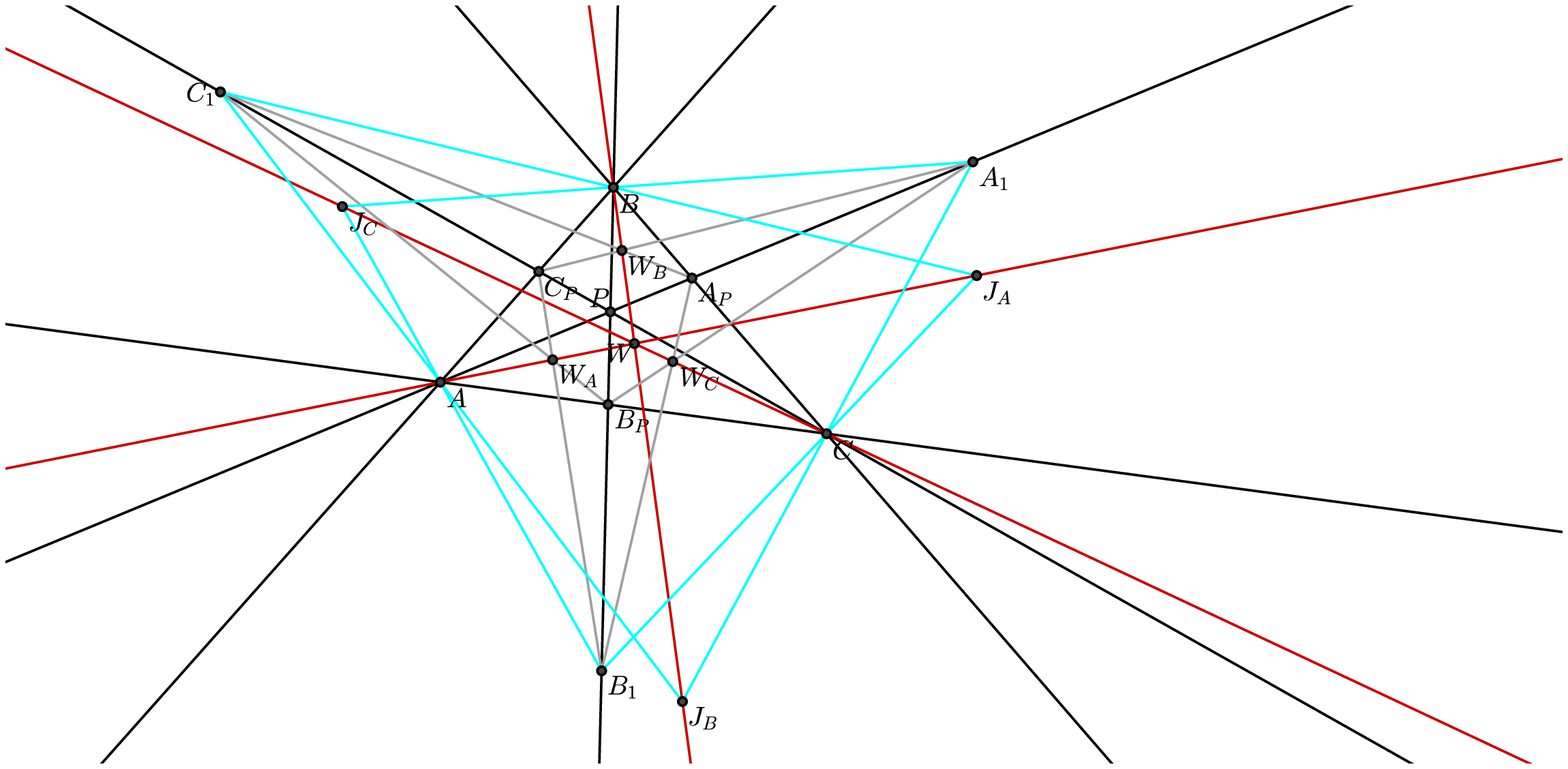}\centering
\caption{Галстук}\label{1}
\end{figure}

\begin{figure}[h]

{\bf Теорема 1.2.} Пусть теперь $A'_P$ --- произвольная точка прямой $AP$. Через $W_B'$ определим точку пересечения $A_1B_P$ и $B_1A'_P$. Аналогично определим $W_C'$. Обозначим через $N$ точку пересечения $BW_B'$ и $CW_C'$. Тогда точки $A$, $B$, $C$, $W$, $A_1$ и $N$ лежат на одной конике.

\smallskip
Доказательство теоремы 1.2 вытекает из следующей леммы.

\smallskip
{\bf Лемма 1.} На плоскости даны пять точек $A$, $B$, $C$, $D$, $E$, и через $E$ проведены 3 прямые $l_1$, $l_2$, $l_3$, не содержащие оставшиеся 4 точки (см. рис. 2). На прямой $l_3$ отметили произвольно точку~$P$. Определим точки $A', B'$ и $X$ так: $A'=CP\bigcap l_1$, $B'=DP\bigcap l_2$, $X=AA'\bigcap BB'$. Тогда точка~$X$ всегда лежит на конкретной конике, проходящей через $A$, $B$ и $E$.

\smallskip
{\it Доказательство леммы 1.}
Определим проективное преобразование $f$ следующим образом. Сначала спроецируем с центром в точке $C$ прямую $l_1$ в прямую $l_2$ и $A$ переведем в себя, а затем спроецируем $l_2$ с центром в $D$ в прямую $l_3$ и оставим опять $A$ на месте, а затем оставим прямую $l_3$ на месте и $A$ переведем в $B$. Очевидно, что так можно сделать, т.\,к. сохраняются двойные отношения и по 4 точкам задается проективное преобразование. Тогда для любой точки $P$ $f(AA')=BB'$ и $f(A)=f(B)$, а значит, в силу леммы Соллертинского (см. [1]),  точка~$X$ принадлежит всегда некоторой конике, которая проходит через $A$ и $B$. Точка $E$ принадлежит этой конике, т.\,к. если $P=E$, то соответствующий $X=E$.
QED

\smallskip
{\it Доказательство теоремы 1.2.} Тогда, очевидно, теорема 1.2 следует из леммы 1. Точкам~$B, C, A_1$ из теоремы $1.2$ соответствуют точки $A$, $B$, $E$ из леммы 1. То, что для точек~$B$ и $W$ найдутся такие положения $A'_P$, что $N$ совпадет с ними следует из теоремы 1.1. QED

\includegraphics{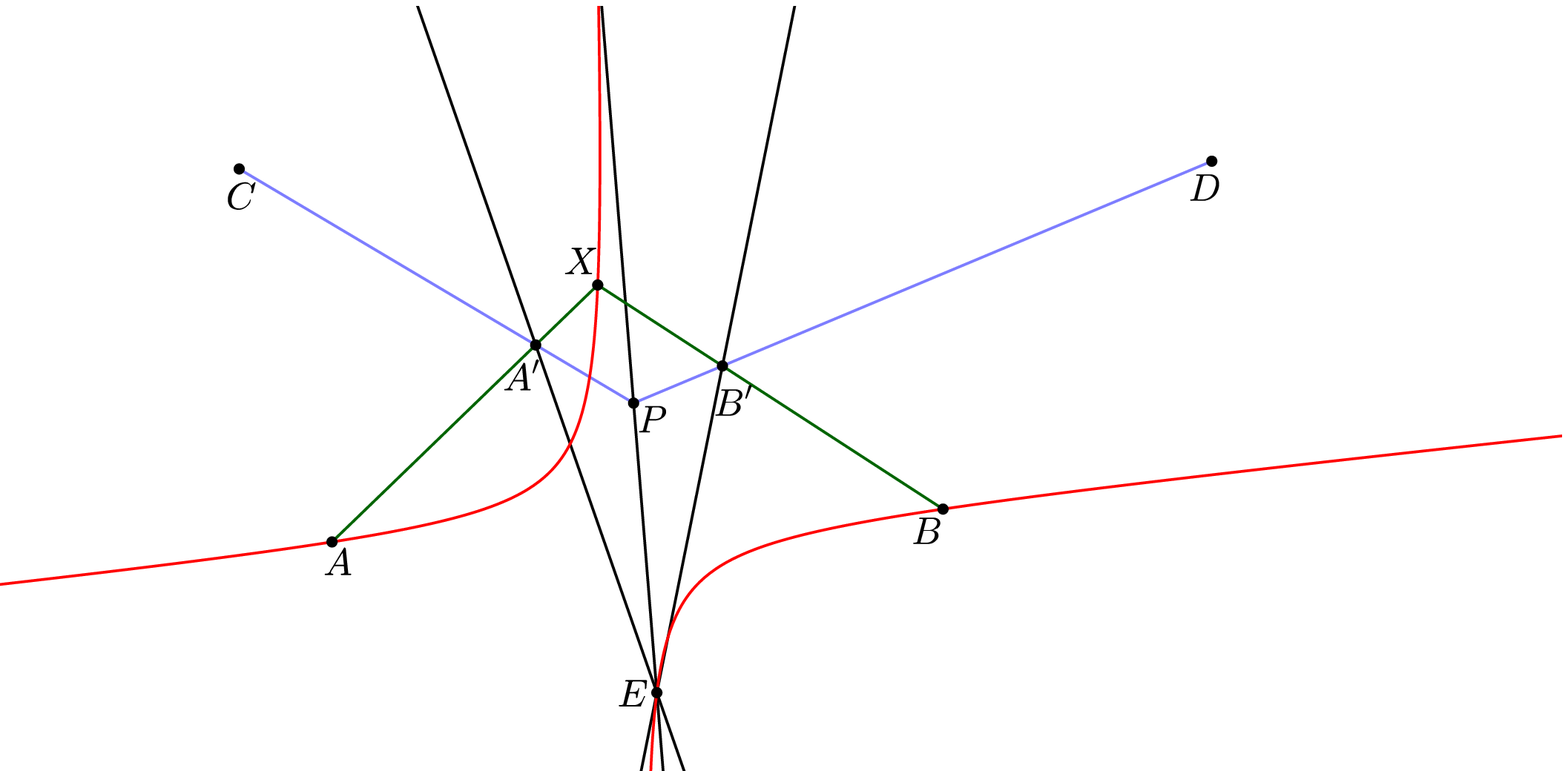}\centering
\caption{Лемма 1}\label{3}
\end{figure}

\begin{figure}[h]
\bigskip
{\bf Теорема 1.3.} Определим через $T_A$ точку пересечения $WA_1$ и $BC$. Аналогично определяются точки $T_B$, $T_C$. Тогда прямые $AT_A$, $BT_B$, $CT_C$ пересекаются в некоторой точке $T$ (см. рис. 3).

\smallskip
{\it Доказательство теоремы $1.3$.} Как мы уже знаем, $C$, $W_C$ и $J_C$ лежат на одной прямой по теореме Паппа. Значит, точки $C$, $W$ и $J_C$ лежат на одной прямой. Тогда по теореме Дезарга треугольники $AT_BB_1$ и $BT_AA_1$ перспективны, так как соответствующие стороны перескаются в точках, лежащих на одной прямой. Это значит, что треугольники $BT_BB_1$ и $AT_AA_1$ перспективны. Значит, $W$, $P$ и точка пересечения $AT_A$ и $BT_B$ лежат на одной прямой. Поскольку аналогичные рассуждения верны для точек пересечения $BT_B$ и $CT_C$, $AT_A$ и $CT_C$, то $AT_A$, $BT_B$, $CT_C$ пересекаются в одной точке. QED

\includegraphics{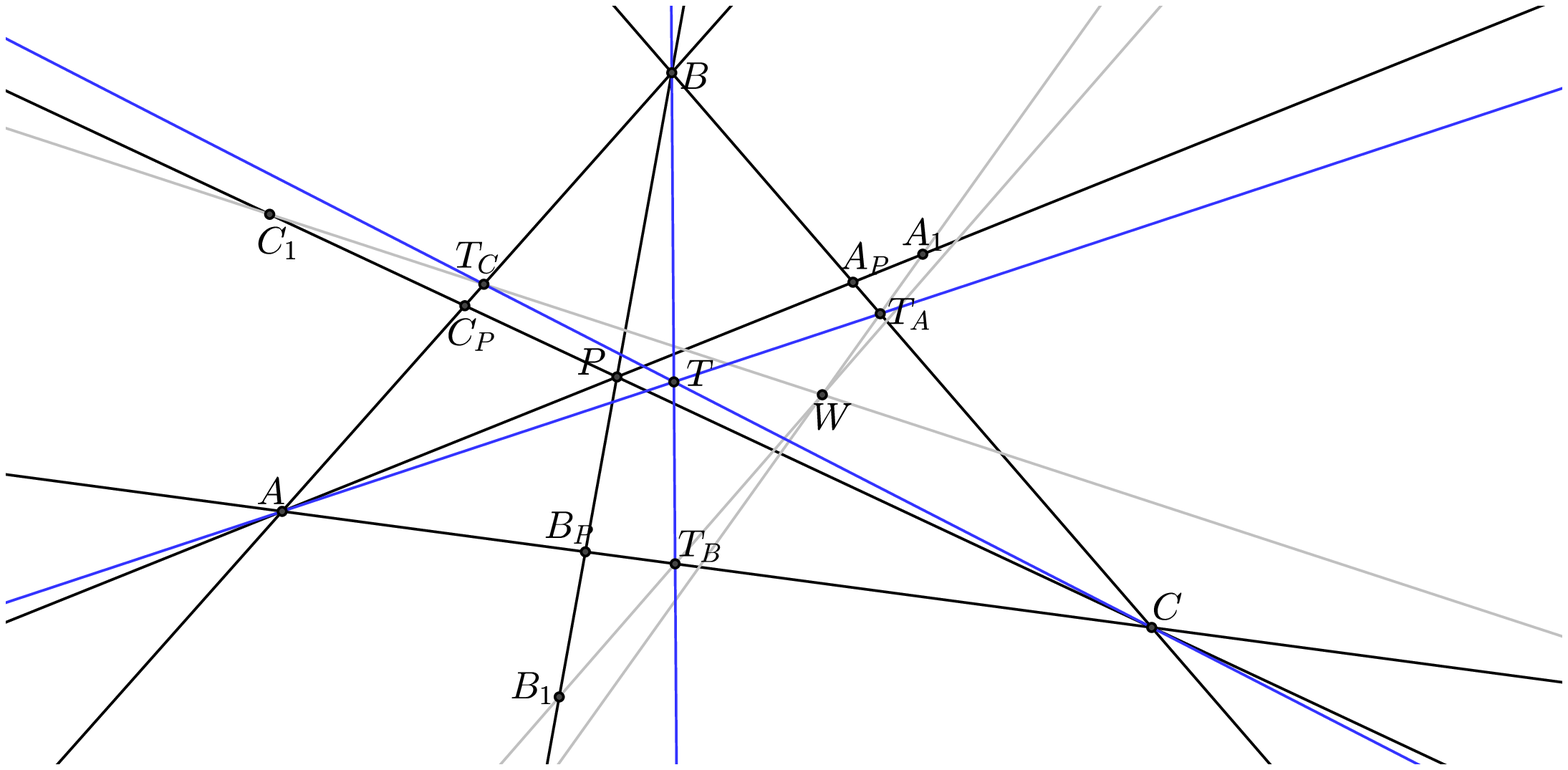}\centering
\caption{Теорема $1.3$}\label{3}
\end{figure}

\bigskip
\begin{figure}[h]
{\bf Теорема 1.4.} За $G_A$ примем точку персечения $W_AP$ и $BC$. Аналогично определяются точки $G_B$ и $G_C$. Тогда $AG_A$, $BG_B$ и $CG_C$ пересекаются в некоторой точке $G$.

\smallskip
{\it Доказательство теоремы $1.4$.} Покажем, что прямые $W_AA_P$, $W_BB_P$ и $W_CC_P$ пересекаются в некоторой точке. В шестиугольник $A_PB_1C_PA_1B_PC_1$ можно вписать конику по обратной теореме Брианшона. Это значит, что эта коника вписана в $A_PW_BC_PW_AB_PW_C$. Значит, по теореме Брианшона указанные выше прямые пересекаются в одной точке. Обозначим эту точку через $R$.

\smallskip
Теперь покажем, что точки пересечения пар прямых $A_PG_B$ и $G_AB_P$, $C_PG_B$ и $G_CB_P$, $A_PG_C$ и $G_AC_P$ лежат на прямой $RP$. Если мы это докажем, то по обратной теореме Паскаля шестиугольник $A_PG_AB_PG_BC_PG_C$ будет вписан в некоторую конику, а значит, $AG_A$, $BG_B$ и $CG_C$ конкурентны. Докажем, что точка пересечения $A_PG_B$ и $G_AB_P$ лежит на $RP$, а для остальных точек доказательство будет аналогичным. Используя теорему Паскаля, получим, что данное утверждение равносильно тому, что $A_PW_AG_AB_PW_BG_B$ вписан в некоторую конику. А это по обратной теореме Паскаля действительно так, потому что $C_1, P, C$ лежат на одной прямой.

{\it Следствие.} Точка $G$ лежит на прямой $RP$.

\includegraphics{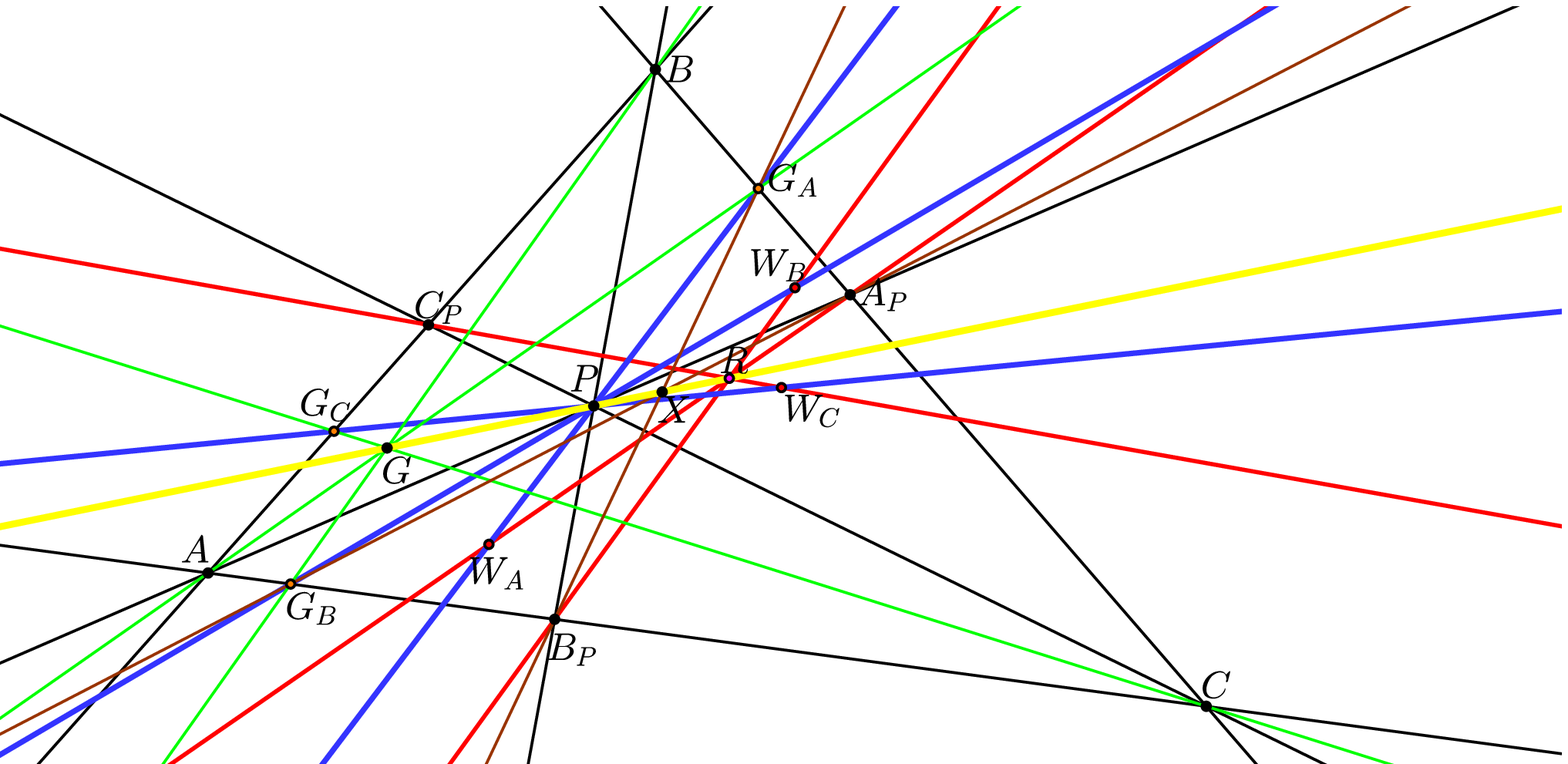}\centering
\caption{Теорема 1.4}\label{4}
\end{figure}

\begin{figure}[h]
{\bf Определения.}

1) ($A, B, C, D$) --- двойное отношение точек $A, B, C, D$.

2) ($l_1, l_2, l_3, l_4$) --- двойное отношение прямых $l_1, l_2, l_3, l_4$.

\smallskip
{\bf Теорема 1.5.} Точки $A$, $B$, $C$, $P$, $G$(из теоремы 1.4) и $T$(из теоремы 1.3) лежат на одной конике.

{\it Доказательство теоремы $1.5$.} Очевидно, что условие теоремы равносильно равенству двойных отношений следующих четверок прямых: $(AB, AG, AP, AT)$ и $(CB, CG, CP, CT)$. Данное равенство, очевидно, эквивалентно тому, что ($B$, $G_A$, $A_P$, $T_A$)=($B$, $G_C$, $C_P$, $T_C$) (так будем обозначать двойное отношение четверки точек). Докажем это равенство.

Прежде заметим, что $PJ_A\bigcap BC=T_A$ (верны и симметричные утверждения для других вершин треугольника). Для этого достатоно заметить, что теорема Дезарга верна для треугольников $PBA_1$ и $J_ACW$ ввиду коллинеарности $A, B_1, J_C$.

 Определим $B_1C\bigcap BA=C_U, B_1A\bigcap BC=A_U, BP\bigcap CW= C_B, BP\bigcap AW= A_B$. Тогда, сделав проекции с центром в точке $P$ с прямых $BC, BA$ на $AW, CW$ соответственно, получим, что ($B$, $G_A$, $A_P$, $T_A$)=($A_B, W_A, A, J_A$) и ($B$, $G_C$, $C_P$, $T_C$)=($C_B, W_C, A, J_C$). Теперь с центром в точке $B_1$ сделаем проекции с прямых $AW$ и $CW$ на $BA$ и $BC$ соответственно. Тогда получим, что ($A_B, W_A, A, J_A$)=($B, C_P, A, C_U$) и ($C_B, W_C, A, J_C$)=($B, A_P, C, A_U$). Теперь заметим, что для треугольников $AA_PA_U$ и $CC_PC_U$ верна теорема Дезарга ввиду коллинеарности $B, P, B_1$. Тогда очевидно, что ($B, C_P, A, C_U$)=($B, A_P, C, A_U$). Значит, ($B$, $G_A$, $A_P$, $T_A$)=($B$, $G_C$, $C_P$, $T_C$). QED

\includegraphics{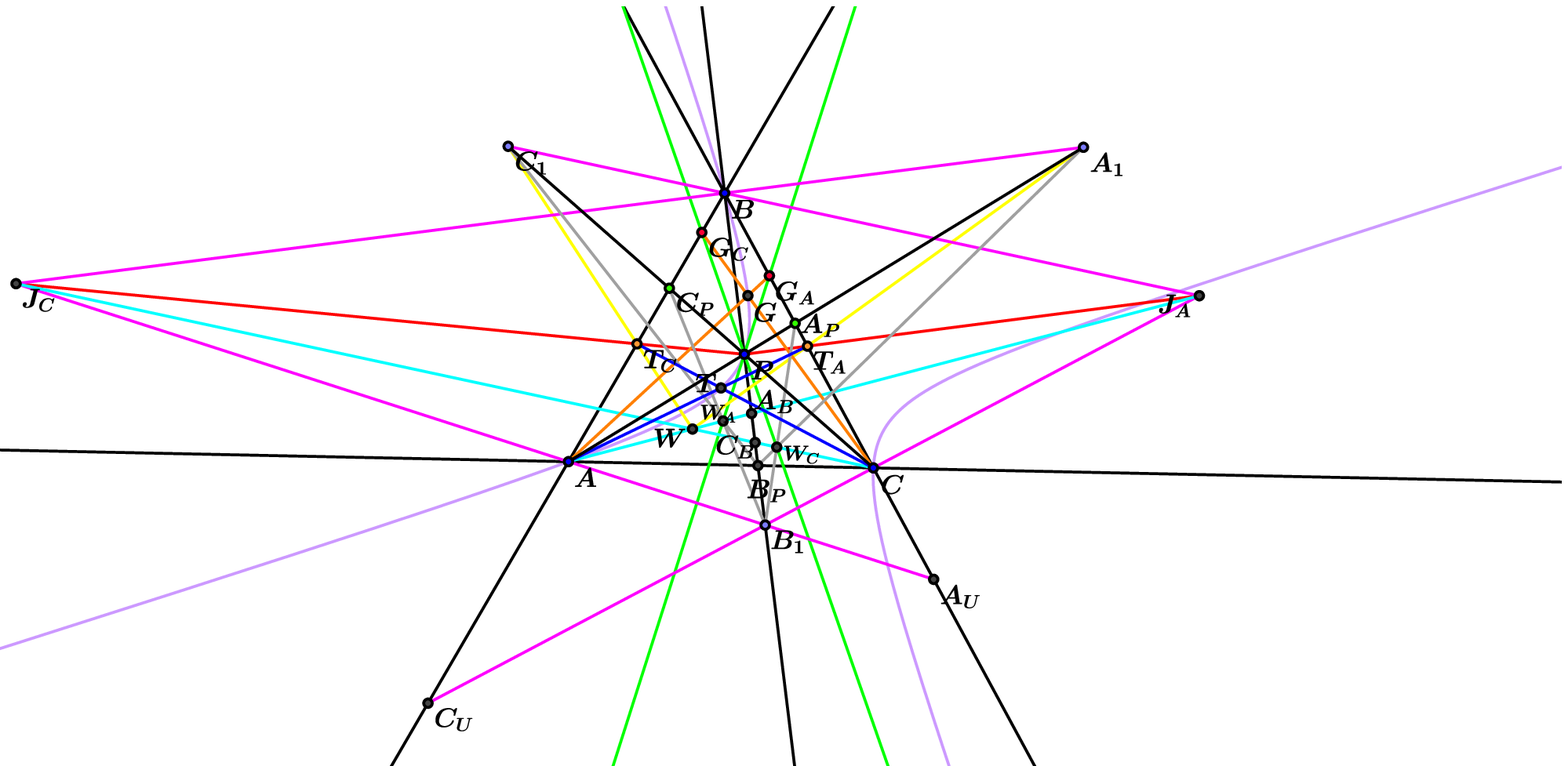}\centering
\caption{Теорема 1.5}\label{5}
\end{figure}\smallskip

\begin{figure}[h]

\bigskip
{\bf Обобщенная теорема о галстуке (Теорема 1.6).} $Q$ --- произвольная точка, не лежащая на сторонах треугольника $ABC$. $A_Q, B_Q, C_Q$ --- точки пересечения прямых $QA_1$, $QB_1$ и $QC_1$ с $BC, AC, AB$ соответственно. $V_C$ --- точка пересечения прямых $A_QB_1$ и $B_QA_1$. Аналогично определяются точки $V_A, V_B$. Тогда $AV_A, BV_B, CV_C$ пересекаются в одной точке.

\smallskip
{\it Доказательство теоремы 1.6.}
$A_1B_QC_1A_QB_1C_Q$ --- описанный шестиугольник по обратной теореме Брианшона. Тогда коника, вписанная в него, вписана и в $A_1V_BC_1V_AB_1V_C$. Следовательно по теореме Брианшона $A_1V_A, B_1V_B, C_1V_C$ пересекаются в некоторой точке. Обозначим ее через $T$.

Заметим, что треугольники $AA_1V_A$ и $BV_BB_1$ перспективны с центром в $C_Q$. Тогда по теореме Дезарга точка $T$ лежит на прямой, проходящей через точки пересечения пар прямых $AV_A$ и $BB_1$, $AA_1$ и $BV_B$, которые определим как $B_A$ и $A_B$ соответственно. Применяя аналогичные рассуждения для пар треугольников $AA_1V_A$ и $CV_CC_1$, $CC_1V_C$ и $BV_BB_1$ и определяя для них точки $A_C, C_A, B_C, C_B$ аналогичным образом, получаем, что $A_BB_A$, $A_CC_A$ и $B_CC_B$ конкурентны.

Теперь вместе с последним результатом из леммы 2 последует сама теорема.

\includegraphics{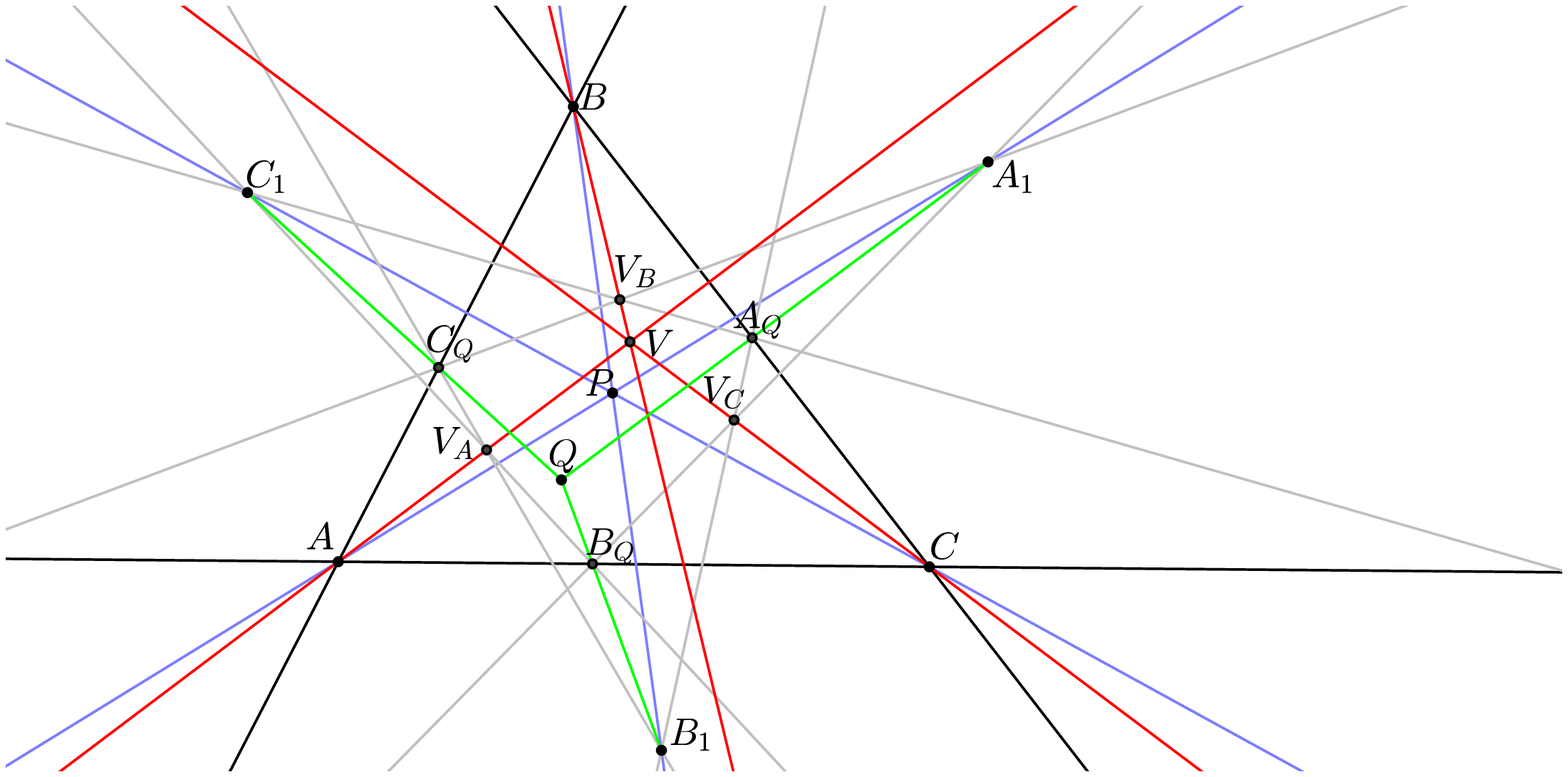}\centering
\caption{Обобщение теоремы о галстуке}\label{6}
\end{figure}

\begin{figure}[h]

  {\bf Лемма 2.} На плоскости есть $5$ различных точек $A, B, C, P, Q$ общего положения. Пусть $AP$ и $BQ$ пересекаются в точке $A_B$, а прямые $AQ$ и $BP$ в точке $B_A$. Аналогично определим точки $B_C, C_B, A_C, C_A$. Тогда прямые $A_BB_A$, $A_CC_A$ и $B_CC_B$ пересекаются в одной точке.

\smallskip
{\it Доказательство леммы 2.}
 Переведем прямую $PQ$ в бесконечно удаленную и афинным преобразованием сделаем направления на $P$ и $Q$ перпендикулярными. Тогда мы получим следующий факт. Имеется прямоугольник $ABCD$. Прямая $l_1$, параллельная $AB$, пересекает $BC$ в точке $K$. Прямая $l_2$, параллельная $AD$, пересекает $AB$ в точке $N$. Прямые $l_1$ и $l_2$ пересекаются в точке $R$. Тогда прямые $NC$, $AK$ и $DR$ конкурентны. Этот факт является общеизвестным и доказывается несложно.
QED

\smallskip
Пусть теперь прямые $AV_A$ и $BV_B$ пересекаются в точке $V$. Мы знаем, что для пары треугольников $B_AC_AC_B$ и $A_BA_CB_C$ верна теорема Дезарга. Тогда точка пересечения $B_AC_B$ и $A_BB_C$ лежит на прямой $AC$. Но $B_C$ как точка, лежащая на прямой $BP$ и удволетворяющая последнему, по точкам $B_A, C_B, A_B$ восстанавливается однозначно. Значит, для пятерки прямых $AP, BP, CP, AV_A BV_B$ существует единственная прямая через точку $C$ , для которой может выполняться конкурентность $A_BB_A$, $A_CC_A$ и $B_CC_B$. По лемме 2 это прямая $CV$. Тогда, очевидно, $CV$ и $CV_C$ совпадают.
QED

\includegraphics{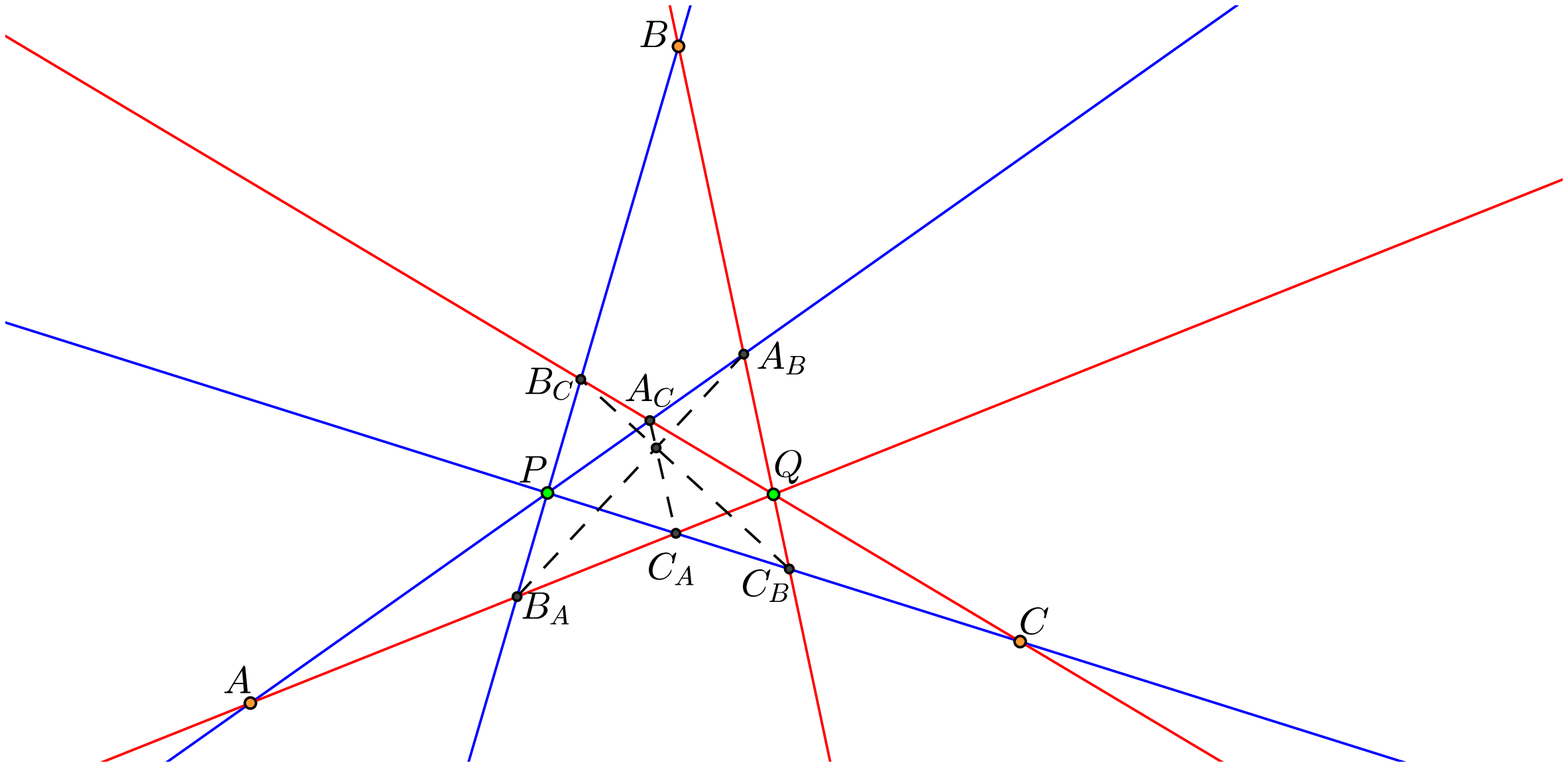}\centering
\caption{Лемма 2}\label{7}
\end{figure}\smallskip

\begin{figure}[h]

\bigskip
{\bf Список литературы.}

\smallskip
1. {\it Акопян А. В., Заславский А. А.} Геометрические свойства кривых второго порядка. 2-е изд., дополн. М.: МЦНМО, 2011

\smallskip
2. {\it А. В. Акопян.} Геометрия в картинках. (c2) М., 2011
\end{figure}

  \end{document}